\documentclass[11pt,reqno]{article}

\usepackage{graphics,color,graphicx}

\usepackage{amssymb}

\usepackage{latexsym}

\usepackage{amsmath}

\newtheorem{theorem}{Theorem}[section]

\newtheorem{question}{Question}[section]

\newtheorem{cor}{Corollary}[section]

\newtheorem{lemma}{Lemma}[section]

\newtheorem{remark}{Remark}[section]

\newtheorem{conj}{Conjecture}[section]

\newtheorem{defin}{Definition}[section]

\newtheorem{claim}{Claim}[section]

\newtheorem{prop}{Proposition}[section]

\newcommand{\proof}{\noindent{\bf\emph{Proof:~~~}}}

\newcommand{\PP}{\mathbb{P}}
\newcommand{\EE}{\mathbb{E}}

\newcommand{\eps}{\varepsilon}

\newcommand{\qed}{\hspace*{\fill}
\rule{7pt}{7pt}}

\DeclareMathOperator {\DEG}{DEG}

\begin{document}

\def\concat#1#2{\sideset{_{#1}}{_{#2}}{\mathop{\circ}}}

\title{Large matchings in uniform hypergraphs
and the conjectures of Erd\H{o}s and Samuels}

\author{
Noga Alon \thanks{Sackler School of Mathematics and Blavatnik School of Computer Science, Tel Aviv
University, Tel Aviv 69978, Israel. Email: {\tt nogaa@tau.ac.il}. Research supported in part by an
ERC advanced grant, by a USA-Israeli BSF grant
and by the Israeli I-Core program.}
\and Peter Frankl \thanks{Tokyo, Japan. Email:
{\tt peter.frankl@gmail.com}} \and Hao Huang
\thanks{Department of Mathematics, UCLA, Los Angeles, CA, 90095. Email: {\tt
huanghao@math.ucla.edu}.} \and Vojtech R\"odl \thanks{Emory University, Atlanta, GA. Email: {\tt
rodl@mathcs.emory.edu}. Research supported  by NSF grant DMS 080070.} \and Andrzej Ruci\'nski
\thanks{A. Mickiewicz University, Pozna\'n, Poland. Email: {\tt rucinski@amu.edu.pl}. Research
supported by the National Science Center grant N N201 604940, and the NSF grant DMS-1102086. Part of research performed at Emory
University, Atlanta.} \and Benny Sudakov\thanks{Department of Mathematics, UCLA,  Los Angeles, CA
90095. Email: {\tt bsudakov@math.ucla.edu}.
Research supported in part by NSF grant DMS-1101185, NSF CAREER award DMS-0812005 and by
USA-Israeli BSF grant.}}

\date{}
\maketitle

\begin{abstract}

In this paper we study degree conditions which guarantee the existence of
perfect matchings and perfect fractional matchings in uniform hypergraphs.
We reduce this problem to an old conjecture by Erd\H{o}s on
estimating the maximum number of edges in a hypergraph when the
(fractional) matching number is given,
which we are able to solve in some special cases using probabilistic techniques.
Based on these results, we obtain
some general theorems on the
minimum $d$-degree ensuring the
existence of perfect (fractional) matchings.
In particular, we asymptotically determine the minimum vertex
degree which guarantees a perfect matching in 4-uniform and
5-uniform hypergraphs.
We also discuss an application to a problem of finding
an optimal data allocation in a distributed storage system.

\end{abstract}

\section{Introduction}
\emph{A $k$-uniform hypergraph} or {\emph a $k$-graph} for short, is a pair $H=(V,E)$, where
$V:=V(H)$ is a finite set of vertices and $E:=E(H)\subseteq  {V \choose k}$ is a family of $k$-element
subsets of $V$ called edges.  Whenever convenient we will identify $H$ with $E(H)$. \emph{A matching} in $H$ is
a set of disjoint edges of~$H$. The number of edges in a matching is called \emph{the size} of the
matching. The size of the largest matching in a $k$-graph $H$ is denoted by $\nu(H)$. A matching is
\emph{perfect} if its size equals $|V|/k$.

\emph{A fractional  matching} in a $k$-graph $H=(V,E)$ is a
function $w:E\to[0,1]$ such that for each $v\in V$ we have
$\sum_{e\ni v}w(e)\le1$. Then $\sum_{e\in E}w(e)$ is the size of
$w$. The size of the largest fractional matching in a $k$-graph
$H$ is denoted by $\nu^*(H)$. If $\nu^*(H)=n/k$, or
equivalently, for all $v\in V$ we have $\sum_{e\ni v}w(e)=1$, then
we call $w$ \emph{perfect}.

The determination of $\nu^*(H)$ is a linear programming problem. Its dual problem is to find a
minimum \emph{fractional vertex cover} $\tau^*(H)=\sum_{v\in V}w(v)$ over all functions
$w:V\to[0,1]$ such that for each $e\in E$ we have $\sum_{v\in e}w(v)\ge1$. Let $\tau(H)$ be the
minimum number of vertices in a vertex cover of $H$. Then, for every $k$-graph $H$, by the Duality
Theorem,
\begin{equation}\label{DT}
\nu(H)\le \nu^*(H)=\tau^*(H)\le\tau(H).
\end{equation}

Given a
$k$-graph $H$ and a set $S\in {V \choose d}$, $0\le d\le k-1$,
we denote by $\deg_H(S)$ the number of edges in $H$ which contain $S$.
 Let $\delta_{d}:=\delta_{d}(H)$ be the minimum $d$-degree of $H$, which is the minimum $\deg_H(S)$ over all $S\in {V \choose d}$.
Note that $\delta_{0}(H)=|E(H)|$. In this paper we study the
relation between the minimum $d$-degree $\delta_{d}(H)$ and the matching numbers $\nu(H)$ and
$\nu^*(H)$.

\begin{defin}\label{m}\rm
  Let integers $d,k,s$, and $n$  satisfy  $0\le d\le k-1$, and $0\le s\le n/k$. We denote by
$m_d^{ s}(k,n)$ the minimum $m$ so that for an $n$-vertex $k$-graph $H$, $\delta_d(H) \geq m$ implies that $\nu(H) \geq s$. Equivalently,
$$m_d^s(k,n)-1=\max\{\delta_d(H): |V(H)|=n\mbox{ and }\nu(H)\le s-1\}.$$
\newline Furthermore, for a \emph{real} number $0\le s\le n/k$, define $f_d^{ s}(k,n)$ as the minimum $m$ so that $\delta_d(H)\ge m$ implies that $\nu^*(H) \ge s$. Equivalently,
$$f_d^s(k,n)-1=\max\{\delta_d(H): |V(H)|=n\mbox{ and }\nu^*(H)<s\}.$$
\end{defin}

Observe that trivially, for $\lceil s\rceil\le n/k$, \begin{equation}\label{flessm}f^s_d(k,n)\le m^{\lceil s\rceil}_d(k,n).\end{equation}

We are mostly interested in the case $s=n/k$ (i.e. when matchings are perfect) in which we suppress the superscript in the notation
$m_d^{n/k}(k,n)$ and $f_d^{n/k}(k,n)$. Thus, writing $m_d(k,n)$, we implicitly require
that $n$ is divisible by $k$.

Problems of this type have a long history going back to Dirac \cite{dirac}
who in 1952 proved that minimum degree $n/2$ implies the existence of a Hamiltonian cycle in graphs. Therefore,
for $d\ge1$, we refer to the extremal parameters $m_d(k,n)$ and $f_d(k,n)$ as to
\emph{Dirac-type thresholds}.
When $k=2$, an easy argument shows that $m_1(2,n)=n/2$. For $k\ge3$, an exact
formula for $m_{k-1}(k,n)$  was obtained in \cite{rrs}.
For a fixed $k\ge3$ and $n\to\infty$ it yields the asymptotics $m_{k-1}(k,n)=\frac{n}{2}+O(1)$.
As far as  perfect fractional matchings are concerned, it was  proved in \cite{RRSm} that $f_{k-1}(k,n)= \lceil n/k\rceil$ for $k\ge2$,
which is a lot less than $m_{k-1}(k,n)$ when $k \geq 3$.
For more results on Dirac-type thresholds
for matchings and Hamilton cycles see \cite{sur}.

In this paper, we focus on the asymptotic behavior of $m_d(k,n)$ and $f_d(k,n)$ for general, but fixed $k$ and $d$, when $n\to\infty$.
For a lower bound on $m_d(k,n)$ consider first a $k$-graph $H_0=H_0(k,n)$ (constructed in
\cite{rrs}) with vertex set split almost evenly, that is, $V(H_0)=A\cup B$, $\big||A|-|B|\big|\le2$, and with the edge set
consisting of all $k$-element subsets of $V(H_0)$ intersecting $A$ in an odd number of vertices. We
choose the size of $A$ so that $|A|$ and $\tfrac nk$ have different parity. Clearly, there is
no perfect matching in $H_0$ and for every $0\le d\le k-1$ we have $\delta_{d}(H_0)\sim \frac{1}{2} {n-d \choose k-d}$.


Another lower bound on $m_d(k,n)$ is given by the following
well known construction.
For integers $n$, $k$, and  $s$, let $H_1(s)$ be a
$k$-graph on $n$ vertices consisting of all $k$-element subsets intersecting a given set of size $s-1$, that is $H_1(s)=K^{(k)}_n-K^{(k)}_{n-s+1}$. Observe that $\nu(H_1(s))= s-1$, while $$\delta_d(H_1(n/k))={n-d \choose k-d}-{n-d-n/k+1 \choose k-d} \sim \left\{1-\left(\frac{k-1}k\right)^{k-d}\right\}{n-d \choose k-d}.$$

Assume that $n$ is divisible by $k$. Putting $s=\tfrac nk$ and using the $k$-graphs $H_0$ and
$H_1(n/k)$, we obtain a lower bound
\begin{equation}\label{genl}
 m_d(k,n)\ge\max\left\{\delta_{d}(H_0),\delta_d(H_1(\tfrac nk))\right\}+1
 \sim\max\left\{\frac12,1-\left(\frac{k-1}k\right)^{k-d}\right\}{n-d \choose k-d}.
\end{equation}
On the other hand, $H_1(\lceil n/k\rceil)$ alone  yields a lower bound also on $f_d(k,n)$. Indeed,  for a real $s>0$ we have
 $$\nu^*(H_1(\lceil s \rceil))=\tau^*(H_1(\lceil s \rceil))\le \tau(H_1(\lceil s \rceil))=\lceil s \rceil -1<s,$$
 and so

\begin{equation}\label{genlf}
f_d(k,n)\ge \delta_d(H_1(\lceil\tfrac nk\rceil))+1
 \sim \left\{1-\left(\frac{k-1}k\right)^{k-d}\right\}{n-d \choose k-d}.
\end{equation}

It is easy to check that for $d\ge k/2$ the  maximum in the coefficient in (\ref{genl}) equals
$\tfrac12$. Pikhurko  \cite{pikh} proved, complementing the case $d=k-1$, that indeed we  have $m_d(k,n)\sim\frac12 {n-d \choose k-d}$ also for $k/2\le d\le k-2$, $k\ge4$.

For $d<k/2$ the problem seems to be harder and we discuss below the cases $d\ge1$ and $d=0$ separately.
The first result for the range $1\le d<k/2$, $k\ge3$, was obtained already in 1981 by Daykin
and H\"aggkvist in \cite{dh} who proved that $m_1(k,n)\le\left ( \frac{k-1}k+o(1)\right ) {n-1 \choose k-1}$.
This was generalized to $ m_d(k,n)\le\left (\frac{k-d}k+o(1)\right ) {n-d \choose k-d}$ for all $1\le d<k/2$ in \cite{hps}, and, using the ideas
from \cite{hps}, slightly improved in \cite{klas} to
$m_d(k,n)\le\left\{\frac{k-d}k-\frac1{k^{k-d}}+o(1)\right\} {n-d \choose k-d}$.
For $k=4, d=1$  the latter coefficient is $\tfrac{47}{64}$. In \cite{klas}, the constant
was further lowered to $\tfrac{42}{64}$, but there is still
a gap between this upper bound and the lower bound of $\tfrac{37}{64}$.

It has been  conjectured in \cite{ko-sur} and
again in \cite{hps} that the lower bound (\ref{genl}) is achieved at least asymptotically.

\begin{conj}\label{M3}
For all $1\le d\le k-1$,
$$m_d(k,n)\sim \max\left\{\frac12,1-\left(\frac{k-1}k\right)^{k-d}\right\}\binom{n-d}{k-d}.$$
\end{conj}

H\`an, Person, and Schacht  in \cite{hps} proved Conjecture \ref{M3} in the case $d=1$,
$k=3$ by showing that
$m_1(3,n)$ is asymptotically equal to $\frac59\binom{n-1}2$.
K\"uhn,
Osthus, and Treglown \cite{kot} and, independently, Khan \cite{khan}, proved the exact result
$m_1(3,n)=\delta_1(H_1(n/3))+1$. Recently Khan \cite{khan4} announced that he verified the exact result $m_1(4,n)=\delta_1(H_1(n/4))+1$, while the asymptotic version, $m_1(4,n)\sim\tfrac{37}{64}\binom{n-1}3$ follows also from a more general result by Lo and Markstr\"om \cite{lo}.

These exact results, together with \eqref{flessm} and \eqref{genlf}, yield that
$f_1(3,n)= m_1(3,n)$ and $f_1(4,n)= m_1(4,n)$.
Remembering that, on the other hand, $f_{k-1}(k,n)$ is much smaller than $m_{k-1}(k,n)$,
one can raise the question about a general relation between $m_d(k,n)$ and its
 fractional counterpart $f_d(k,n).$ In this paper we
answer this question
 by showing that $m_d(k,n)$ and  $f_d(k,n)$ are  asymptotically equal whenever $f_d(k,n)\sim c^*\binom{n-d}{k-d}$
 for some constant $c^*>\tfrac12$, and otherwise $m_d(k,n)\sim
 \tfrac12\binom{n-d}{k-d}$.

\begin{theorem}\label{rel} For every $1\le d\le k-1$ if there exists $c^*>0$ such that
 $f_d(k,n)\sim c^*\binom{n-d}{k-d}$ then
\begin{equation}\label{sim}
m_d(k,n)\sim\max \left \{c^*,\tfrac12 \right\}\binom{n-d}{k-d}.
\end{equation}
\end{theorem}

This result reduces the task of asymptotically calculating $m_d(k,n)$ to a presumably simpler task
of calculating $f_d(k,n)$.
It seems that, similarly to the integral case, the lower bound
in (\ref{genlf}) determines asymptotically the actual value of the parameter $f_d(k,n)$.
\begin{conj}\label{F3}
For all $1\le d\le k-1$,
$$f_d(k,n)\sim \left\{1-\left(\frac{k-1}k\right)^{k-d}\right\}\binom{n-d}{k-d}.$$
\end{conj}
Our next result confirms Conjecture \ref{F3}  asymptotically for all $k$ and $d$ such that $1\le
k-d\le4$.
Note that the above mentioned result from \cite{RRSm} shows that Conjecture \ref{F3} is true for $d=k-1$ exactly, that is, $f_{k-1}(k,n)=\delta_{k-1}\left(H_1\left(\lceil\tfrac nk\rceil\right)\right)+1$. We include this case into the statement of Theorem \ref{exact} for completeness.

\begin{theorem}\label{exact} For every $k\ge3$ and $k-4\le d\le  k-1$, we have
$$f_{d}(k,n)\sim \left\{1-\left(\frac{k-1}k\right)^{k-d}\right\}\binom{n-d}{k-d}.$$
\end{theorem}

Theorems \ref{exact} and \ref{rel} together imply immediately the
validity of Conjecture \ref{M3} in  a couple of new instances (as discussed earlier, the first of them has been recently also proved in \cite{khan4} and \cite{lo}).

\begin{cor}\label{4152} We have
\begin{enumerate}
\item[] $m_1(4,n)\sim \tfrac{37}{64}\binom{n-1}3$, \quad $m_2(5,n)\sim \tfrac12\binom{n-2}3$, \quad $m_1(5,n)\sim \tfrac{369}{625}\binom{n-1}4$

\item[] $m_2(6,n)\sim \tfrac{671}{1296}\binom{n-2}4$, \quad $m_3(7,n)\sim \tfrac12\binom{n-3}4$.

\end{enumerate}
\end{cor}

We prove Theorem \ref{exact} utilizing the following connection between the parameters $f_d^s(k,n)$ and $f_0^s(k-d,n-d)$.

\begin{prop}\label{ff}
For all $k\ge3$, $1\le d\le k-1$, and $n\ge k$,
$$ f_d(k,n)\le f_0^{n/k}(k-d,n-d).$$
\end{prop}
In view of Proposition \ref{ff}, in order to prove Theorem \ref{exact} we need to estimate
  $f^s_0(k-d,n-d)$ with  $s=\tfrac nk$.
This is trivial for $d=k-1$ and so, from now on, we will be assuming that $d\le k-2$.
The integral version  of this problem has almost as long history as the Dirac-type problem ($d\ge1$).

Erd\H os and Gallai \cite{EG} determined $m^s_0(k,n)$ for graphs ($k=2$). In 1965, Erd\H os
\cite{E}  conjectured the following hypergraph generalization of their result.

\begin{conj}\label{e}
For all $k \ge 2$ and $1\le s\le \tfrac nk$:
\begin{equation*}
m^s_0(k,n)=\max\left\{ \binom{ks-1}k, \binom nk-\binom{n-s+1}k \right\}+1.
\end{equation*}
\end{conj}

The lower bound comes from considering again the extremal $k$-graph $H_1(s)$ along with the
$k$-uniform clique $K_{ks-1}^{(k)}$ (complemented by $n-ks+1$ isolated vertices) which, clearly, has no matching of size $s$. For more on Erd\H os'
conjecture we refer the reader to the survey paper \cite{F-sur} and a recent paper \cite{3}, where
the conjecture is proved for $k=3$ and $n\ge 4s$. In its full generality, the conjecture is still wide open.

We now formulate the  fractional version of Erd\H os' Conjecture.  For future references, we switch
from $k$ and $n$ to $l$ and $m$. Again, the lower bound is yielded by  $H_1(\lceil
s\rceil)$ and the complete $l$-graph on $\lceil ls\rceil-1$ vertices, $K_{\lceil ls\rceil-1}^{(l)}$.

\begin{conj}\label{Efrac}
For all integers $l\ge2$ and an integer $s$ such that  $0\le s\le
m/l$, we have
$$f^s_0(l,m)=\max\left\{ \binom{\lceil ls\rceil-1}l, \binom ml-\binom{m-\lceil s\rceil+1}l \right\}+1.$$
\end{conj}

Note that Conjecture \ref{Efrac} implies that the bound is also asymptotically true for non-integer values of $s$, when $m$ is large.
In \cite{dstorage_arxiv}, there is an example showing that the stronger, precise
version of the conjecture does not hold for fractional $s$.

As a consequence  of the Erd\H os-Gallai theorem from \cite{EG},  Conjecture \ref{Efrac} is asymptotically true for $l=2$ and $m$ goes to infinity.
In the next section we establish a result which confirms Conjecture \ref{Efrac}
asymptotically in the two smallest new instances, but limited to the range $0 \leq s \leq \tfrac m{l+1}$. In this range the case $l=3$ follows
also from the above mentioned result in \cite{3}. It is easy to check that for $s\le \tfrac m{l+1}+O(1)$, the maximum in Conjecture \ref{Efrac} is achieved by the second term.

\begin{theorem}\label{34} For  $l\in\{3,4\}$, for all $d\ge 1$, and $s=\tfrac{m+d}{l+d},$
$$f^{s}_{0}(l,m)\sim \left\{1-\left(1-\frac{1}{l+d}\right)^{l}\right\}\binom{m}{l}$$ where the
asymptotics holds for $m\to\infty$ with $d$ fixed.
\end{theorem}

Theorem \ref{34} together with Proposition \ref{ff} implies Theorem \ref{exact}, which, in turn, together with Theorem \ref{rel} yields Corollary \ref{4152}.
To prove Conjecture \ref{M3} in full generality, one would need to prove Theorem \ref{34} for all $l$.

The rest of this paper is organized as follows.
In the next section, we prove
Theorem \ref{34} using as a main tool a probabilistic inequality of Samuels.
A proof of Proposition \ref{ff}, and consequently of Theorem \ref{exact},  appears in Section \ref{proofexact}. Section \ref{proofrel}  contains a proof of
Theorem \ref{rel}. Finally,  in Section \ref{dstorage}, we  discuss an application of the fractional version of the Erd\"{o}s problem in distributed storage allocation. The last section contains concluding remarks and open problems.

\section{Fractional matchings and
probability of small deviations}\label{proof34}

In this section
we prove Theorem \ref{34}
using  a probabilistic approach from \cite{ahs} based on a special
case of an old probabilistic conjecture of Samuels
\cite{samuels_chebyshev}.
In fact, we prove a little bit more -- see Corollary \ref{t95} and
Remark \ref{r21} below.


For $l$  reals $\mu_1,\ldots,\mu_l$
satisfying $0\le\mu_1 \leq \mu_2 \leq \cdots \leq \mu_l$ and
$\sum_{i=1}^l \mu_i <1$, let
$$P(\mu_1,\mu_2, \ldots ,\mu_l)=\inf\PP(X_1+\ldots +X_l<1),$$ where the infimum is taken
over all possible collections of $l$ independent nonnegative random
variables $X_1, \ldots ,X_l$, with expectations $\mu_1, \ldots
,\mu_l$, respectively. Define
$$Q_t(\mu_1, \ldots ,\mu_l)
=\prod_{i={t+1}}^l\left(1-\frac{\mu_i}{1-\sum_{j=1}^t \mu_j}\right)$$ for each $0 \leq t < l.$

Note that $Q_t(\mu_1, \ldots ,\mu_l)$ is exactly  $\PP(X_1+\ldots +X_l<1)$ when
$X_i$ is identically $\mu_i$ for all $i \leq t$, while $X_i$ attains the values $0$ and
$1-\sum_{i\leq t}\mu_i$ (with its expectation being $\mu_i$) for all $i\ge t+1$.

The following conjecture was raised by Samuels in
\cite{samuels_chebyshev}.
\begin{conj}[\cite{samuels_chebyshev}]
\label{c92}
For all admissible values of $\mu_1, \ldots ,\mu_l$,
$$
P(\mu_1,\mu_2, \ldots ,\mu_l)=\min_{t=0,\dots,l-1} Q_t(\mu_1,\mu_2,
\ldots ,\mu_l).
$$
\end{conj}
Note that for $l=1$ this is Markov's inequality.
Samuels proved his conjecture for $l \leq 4$ in \cite{samuels_chebyshev, samuels_more}.
\begin{lemma}[\cite{samuels_chebyshev, samuels_more}]
\label{l93}
The assertion of Conjecture \ref{c92} holds for all $l \leq 4$.
\end{lemma}
We next show that for $\mu_1=\mu_2= \cdots =\mu_l=x$, where $0 < x \leq \frac{1}{l+1}$, the minimum in Conjecture \ref{c92} is attained by $Q_0(\mu_1, \ldots ,\mu_l).$
\begin{prop}
\label{l94}
For every integer $l \geq 2$ and every  real number $x$ satisfying
$0 < x \leq \frac{1}{l+1}$, if $\mu_1=\mu_2= \ldots
=\mu_l=x$ then
$$
\min_{t=0,\dots,l-1} Q_t(\mu_1,\mu_2,
\ldots ,\mu_l)=Q_0(\mu_1,\mu_2, \ldots ,\mu_l)=(1-x)^l.
$$
\end{prop}
\proof
By definition
$$
Q_t(\mu_1,\mu_2, \ldots ,\mu_l)=\Big(1-\frac{x}{1-tx}\Big)^{l-t}
=\Big(\frac{1-(t+1)x}{1-tx}\Big)^{l-t}.
$$
We thus have to prove that for $0 < x \leq \frac{1}{l+1}$ and $1 \leq t \leq l-1$,
$$
(1-x)^l \leq \Big(\frac{1-(t+1)x}{1-tx}\Big)^{l-t}
$$
or equivalently that
$$
\Big(\frac{1}{1-x}\Big)^l \geq \Big(\frac{1-tx}{1-(t+1)x}\Big)^{l-t}.
$$
By the geometric-arithmetic means inequality applied to
a set of $l$ numbers, $t$ of which are equal to $1$ and the
remaining $l-t$ equal to the quantity
$\frac{1-tx}{1-(t+1)x}$,
we conclude that
$$
\Big(\frac{1-tx}{1-(t+1)x}\Big)^{l-t} \cdot 1^t
\leq
\Big[\frac{1}{l} \cdot \Big(\frac{(1-tx)(l-t)}{1-(t+1)x} +t \Big)\Big]^l.
$$
Thus, it suffices  to show that
$$
\frac{(1-tx)(l-t)}{1-(t+1)x} +t \leq \frac{l}{1-x}.
$$
This is equivalent to
$$
(1-x)[(1-tx)(l-t)+t-t(t+1)x] \leq l[1-(t+1)x],
$$
which is equivalent to
$$
(1-x)[l-t(l+1)x] \leq l-l(t+1)x,
$$
or
$$
l-t(l+1)x -lx +t(l+1)x^2 \leq l-l(t+1)x.
$$
After dividing by $x$, we see that this is equivalent to $x \leq \frac{1}{l+1}$,
which holds by assumption, completing the proof. \qed


Note that when $s=xm$ and $x \leq \frac{1}{l+1}$,
the maximum in Conjecture
\ref{Efrac} is achieved by the second term. We now prove the following, in most part conditional, result, which
shows how to deduce Conjecture \ref{Efrac} in this range from Conjecture \ref{c92}.

\begin{theorem}\label{t94}
For any $l \geq 3$ and $0< x \leq \frac{1}{l+1}$, if Conjecture \ref{c92} holds for  $l$
and $\mu_1=\mu_2 = \ldots =\mu_l=x $ then

$$
f_0^{xm}(l,m) \sim \left\{1-(1-x)^l\right\} {m \choose l}.
$$
\end{theorem}

Combining Theorem \ref{t94} with Lemma \ref{l93}, we obtain the following corollary which implies  Theorem \ref{34}. (For $d=1$, observe that $f^s_0(l,m)\sim f^s_0(l,m+1)$.)

\begin{cor} \label{t95}
For $l=3$, $x \leq 1/4$ and for $l=4$, $x \leq 1/5$,
the maximum number of edges in an $l$-uniform hypergraph $H$
on $m$ vertices with fractional matching number less than $xm$ is $$f_0^{xm}(l,m) \sim\left\{1-(1-x)^l\right\}{m \choose l}.$$
\end{cor}

\textbf{\emph{Proof of Theorem \ref{t94}}}:
Let $H$ be an $l$-uniform hypergraph on a vertex set $V$, $|V|=m$, and suppose
that $\nu^{*}(H)<xm$.
By duality, we also have $\tau^{*}(H)<xm$,
and hence there exists a  weight function $w:V\to [0,1]$
such that $\sum_{v\in V} w(v)<xm$ and, for
every edge $e$ of $H$, $\sum_{v\in e} w(v) \geq 1$. By increasing
the weights $w(v)$ if needed, we may assume that
$$\sum_{v\in V}^mw(v)=xm.$$
  Let $v_1,\dots, v_l$ be a sequence of random vertices of $H$, chosen independently and uniformly at random from $V$.
For each $i=1,\dots,l$ we define a random variable $X_i=w(v_i)$.
Note that $X_1, X_2, \ldots ,X_l$ are independent,
identically distributed random variables, where every
$X_i$ attains each of the $m$ values $w(v)$ with probability $1/m$.
(The values of $w$ for different vertices can be equal, but this is of no importance for us.)

By definition, the expectation $\mu_i$ of each $X_i$ is
$$\mu_i=\sum_{v\in V}\frac1m\cdot w(v)=\frac{xm}m=x.$$

Now we can estimate the number of edges of $H$ as follows.
Since for each edge of $H$ we have $\sum_{v\in e} w(v) \geq 1$, the number $N$ of all $l$-element subsets $S$ of $V$ with
$\sum_{v\in S} w(v) < 1$ is a lower bound on the number of non-edges of $H$.
Let $N_1$ and $N_2$ be the numbers of  all $l$-element sequences of vertices of $V$ and all  $l$-element sequences of \emph{distinct} vertices of $V$, respectively, with the sums of weights strictly smaller than 1.
Note that $N_1-N_2$ is at most the number of $l$-element sequences in which at least one vertex appears twice, thus it is bounded by $\binom{l}{2} m^{l-1} = O(m^{l-1}).$
As the number of all $l$-element subsets of $V$ is $\binom ml=(1+o(1))m^l/l!$ and $N=N_2/l!$,
we have
$$\PP\left(\sum_{i=1}^lw(v_i)<1\right)=\frac{N_1}{m^l}\le\frac{N_2+O(m^{l-1})}{\binom mll!}=(1+o(1))\frac{N}{\binom ml}.$$
If Conjecture \ref{c92} holds for a given $l$ then,
by Lemma \ref{l93} and Proposition \ref{l94},
$$\PP\left(\sum_{i=1}^lw(v_i)<1\right)=\PP\left(\sum_{i=1}^lX_i<1\right)\ge(1-x)^l,$$
and, consequently,
$$N\ge(1+o(1))(1-x)^l\binom ml.$$
It follows that the number of
edges of $H$ is at most
$(1+o(1))\left\{1-(1-x)^l\right\} {m \choose l}$, as needed.
\qed



\begin{remark}
\label{r21}
\rm Note that the above proof works
as long as the conclusion of Proposition \ref{l94} holds.  One
can check  using Mathematica that Proposition \ref{l94} holds for $l=3$ and all $0 < x \leq 0.277$, as well as for $l=4$ and all  $0 < x  \leq 0.217$. Therefore, Corollary \ref{t95} extends to these broader ranges of $x$.
 For bigger values of $x$, e.g., for $x=0.3$
when $l=3$, this is not the case anymore, and the above method
does not suffice to determine the asymptotic behavior of
$f_0^{xm}(l,m)$. In fact, using Samuels conjecture
in the higher range of $x$, one gets a bound on $f_0^{xm}(l,m)$ which is larger than that in Conjecture \ref{Efrac}.
 However, in view of Proposition \ref{ff}, for our main application the case $x \leq
\frac{1}{l+1}$  is just what we need.

\end{remark}

\section{Thresholds for perfect fractional matchings}\label{proofexact}

In this section we present a proof of Proposition \ref{ff} and then deduce quickly Theorem \ref{exact}.

\textbf{\emph{Proof of Proposition \ref{ff}}}:
The outline of the proof goes as follows. We will assume that there is no fractional
perfect  matching in a $k$-graph $H$ on $n$ vertices and then show that the neighborhood graph $H(L)$  in $H$ of a particular set $L$ of size $d$ satisfies $\nu^*(H(L))<n/k$. This will imply that
 $\delta_d(H)\le|H(L)|<f_0^{n/k}(k-d,n-d)$. In contrapositive,
we will prove that if $\delta_d(H)\ge f_0^{n/k}(k-d,n-d)$ then $H$ has a fractional perfect matching, from which it follows, by definition, that $f_d(k,n)\le f_0^{n/k}(k-d,n-d)$.

Let an $n$-vertex $k$-graph $H$ satisfy $\nu^*(H)<n/k$, that is, have no fractional perfect matching. As $\tau^*(H)=\nu^*(H)$, there is a function $w: V\to[0,1]$ such that $\sum_{v\in V}w(v)<n/k$ and, for every $e\in H$, we have $\sum_{v\in e}w(v)\ge1$.
 We can replace $H$ with the $k$-graph whose edge set consists of every $k$-tuple of vertices
on which $w$ totals to at least one.

Formally,  for every weight function $w:V\to[0,1]$ define
$$H_{w} :=\left\{e\in\binom V{k}: \sum_{v\in e}w(v)\ge 1\right\}.$$

For a given weight function $w$, suppose $L$ is a set of $d$ vertices with the smallest weights.
Without loss of generality, we may assume that the $d$ lowest values of $w(x)$ are all
equal to each other, since otherwise we could replace them by their average. (Obviously, this modification would not
change $\sum_{v \in V} w(v)$ nor the set $L$.)
Note that the minimum $d$-degree $\delta_d(H_w)=\min_{S \subset \binom{V}{d}} \deg_H(S)$ is achieved
when $S=L$. Let $H(L)$ be the neighborhood of $L$ in $H_w$, that is a $(k-d)$-graph on the vertex set $V\setminus L$ and with the edge set
$$\left \{S \in \binom{V- L}{k-d}: S \cup L \in E(H_w)\right\}.$$
Then $|H(L)|=\delta_d(H_w)$ and it remains to prove that $\tau^*(H(L))<n/k.$

Let $w_0=\min_{v\in V}w(v)$ and observe that $w_0<1/k$. If $w_0>0$, apply to the weight function $w$ the following linear map
$$w'=\frac{w-w_0}{1-kw_0}.$$
Then, still $\sum_{v\in V}w'(v)<n/k$ and $H_w=H_{w'}$. Moreover, for every $v\in L$, we have $w'(v)=0$.
It follows that the function $w'$ restricted to the set
$V\setminus L$ is a fractional vertex cover of $H(L)$ and
so $\nu^*(H(L))=\tau^*(H(L))<n/k$,
which completes the proof of Proposition \ref{ff}. \qed

\textbf{\emph{Proof of Theorem \ref{exact}}}:
As explained earlier, $f_0^{n/k} (k-d,n-d)=n/k$  holds trivially for $d=k-1$ and together with Proposition  \ref{ff} implies the theorem in this case.
For $d=k-2$, we apply Proposition  \ref{ff} together with the case $l=2$ of the  fractional Erd\H{o}s Conjecture \ref{Efrac} (as mentioned earlier, it follows asymptotically from \cite{EG}).
For $d=k-3$ and $d=k-4$, we use Proposition \ref{ff} and Corollary \ref{34} proved in the previous section. \qed

\begin{remark}\rm Consider a restricted version of Samuels' problem to minimize $\PP(X_1+\cdots +X_l<1)$ under the \emph{additional} assumption that all random variables are identically distributed. Our proofs
indicate that under this regime, for a given $l\ge5$ and $\mu_1=\cdots=m_l=x\le\tfrac1{l+1}$, if
$$\PP(X_1+\cdots +X_l<1)\ge (1+o(1))(1-x)^l$$
then Theorem \ref{exact} would hold for all $k\ge l+1$ and $d=k-l$.

\end{remark}

\section{Constructing integer matchings from fractional ones}\label{proofrel}
In this section, we will prove Theorem \ref{rel}. An indispensable tool in our proof is the Strong Absorbing
Lemma~\ref{SA} from \cite{hps} (see Lemma 10 therein). This lemma provides a sufficient condition on degrees and co-degrees of a hypergraph ensuring the existence
of a small and powerful matching which, by ``absorbing'' vertices, creates a perfect matching from any
 nearly perfect matching.

\begin{lemma}\label{SA}
 For all $\gamma > 0$ and integers $k > d> 0$ there is an $n_0$
such that for all $n > n_0$ the following holds: suppose that $H$ is a $k$-graph
on $n$ vertices with $\delta_d(H)\ge (1/2 + 2\gamma)\binom{n-d}{k-d}$, then there exists a
matching $M:=M_{abs}$ in $H$  such that
\begin{enumerate}
\item[(i)] $|M|<\gamma^kn/k$, and
\item[(ii)] for every set $W \subset V\setminus V (M)$ of
size at most $|W|\le \gamma^{2k}n$ and divisible by $k$ there exists a matching in $H$ covering exactly the vertices
of $V (M) \cup W$.
\end{enumerate}
\end{lemma}

Equipped with this lemma we can practically reduce our task to finding an almost perfect matching
in a suitable subhypergraph of $H$. Here is an outline of our proof of Theorem~\ref{rel}. Assume that
there exists a constant $0< c^*<1$ such that $f_{d}(k,n)\sim c^*\binom{n-d}{k-d}$.
 For any $\alpha>0$ consider a $k$-graph $H$ on $n$ vertices, where $n$ is sufficiently large, with
$$\delta_d(H)\ge(c+\alpha)\binom{n-d}{k-d},$$
where $c=\max \{\tfrac12,c^*\}$. Our goal is to show that $H$ contains a perfect matching.

Set  $\gamma =\alpha/2$ and $\eps=\gamma^{2k}$.
The proof consists of three steps.

\begin{enumerate}
\item Find an absorbing matching $M_{abs}$ satisfying properties (i) and (ii) of  Lemma \ref{SA}.
  Set $H'=H\setminus V(M_{abs})$ and note that when $n$ is sufficiently large,
  $$\delta_d(H') \ge \delta_d(H)-\left(\binom{n-d}{k-d}-\binom{n-d-\varepsilon n}{k-d}\right)
  \ge(c+\alpha/2)\binom{n-d}{k-d} = (c+\gamma)\binom{n-d}{k-d}.$$

\item Find a matching $M_{alm}$ in $H'$ such that $|V(M_{alm})|\ge (1-\eps)|V(H')|$, and thus,
$|V(M_{alm}\cup M_{abs})|\ge (1-\eps)n$.
\item Extend $M_{alm}\cup M_{abs}$  to a perfect matching of $H$ by using
 the absorbing property (ii) of $M_{abs}$ with respect to $W=V(H')\setminus V(M_{alm})$.
\end{enumerate}

Now come the details of the proof.
The Strong Absorbing Lemma provides an absorbing matching $M_{abs}$, so Steps 1 and 3 are clear.
Hence to complete the proof of Theorem \ref{rel} it remains to explain Step 2. One possible approach to
find an almost perfect matching in $H'$ is via the weak hypergraph regularity lemma. Our proof,
however, is based on Theorem~1.1 in \cite{fr}. Recall that the
$2$-degree of a pair of vertices in a hypergraph is the number of
edges containing this pair.
An immediate
corollary of  that theorem  asserts the
existence of an almost perfect matching in any nearly regular $k$-graph
in which all $2$-degrees
are much smaller than
the vertex degrees.
(See Remark after Theorem 1.1 in \cite{fr} or Chapter 4.7 of \cite{AS}).
Here we
formulate this corollary as the following lemma in which $\Delta_2(H)$ denotes the maximum $2$-degree in $H$.

\begin{lemma}\label{FR}
For every integer $k \geq 2$ and a real $\varepsilon>0$, there exists $\tau=\tau(k,\varepsilon)$, $d_0=d_0(k,\varepsilon)$ such that
for every $n  \geq D \geq d_0$ the following holds.

Every $k$-uniform hypergraph on a set $V$ of $n$ vertices which satisfies the following conditions:
\begin{enumerate}
\item $(1-\tau)D< \deg_H(v)<(1+\tau)D$ for  all $v\in V$, and
\item $\Delta_2(H)<\tau D$
\end{enumerate}
contains a matching $M_{alm}$ covering all but at most $\eps n$ vertices.
\end{lemma}

Hence, Step 2 above reduces to finding a spanning subhypergraph $H''$ of $H'$  satisfying the
assumptions of Lemma \ref{FR} with $\eps=\gamma^{2k}$ and
 other parameters $\tau, D, a$ to be suitably chosen. Indeed, the following claim is all we need to complete the proof of Theorem \ref{rel}.
For convenience, we set $n:=|V(H')|$. Recall that $c= \max\{\frac{1}{2}, c^*\}$ where $c^{*}$ comes from the threshold which guarantees the existence of fractional perfect matchings.

\begin{claim}\label{all} For sufficiently large $n$, any $k$-graph $H'$ on $n$ vertices satisfying
$\delta_d(H')\ge(c+\gamma)\binom{n-d}{k-d}$ contains a spanning subhypergraph $H''$, such that for all $v\in V(H'')$ we have $\deg_{H''}(v)\sim n^{0.2}$ while $\Delta_2(H'')\le n^{0.1}$.
\end{claim}

Consequently for every $k\ge2$, $\eps>0$, the subhypergraph $H''$  satisfies the assumptions of Lemma \ref{FR} with $D=n^{0.2}$, and any $\tau>0$.
We obtained the following result as an immediate corollary, which asserts the validity of Step 2 and completes our proof of Theorem \ref{rel}.

\begin{cor}
$H'$ contains an almost perfect matching covering at least $(1-\varepsilon)|V(H')|$ vertices.
\end{cor}

In the proof of Claim \ref{all},
the following well-known concentration results (see, for example \cite{AS}, Appendix A,  and Theorem 2.8, inequality (2.9) and (2.11) in
\cite{JLR}) will be used several times.
We denote by $Bi(n,p)$ a binomial random variable with parameters $n$ and $p$.
\begin{lemma} (Chernoff Inequality for small deviation) If $X = \sum_{i=1}^n X_i$, each random variable
$X_i$ has Bernoulli distribution with expectation $p_i$, and $\alpha \leq 3/2$, then
\begin{equation}\label{smalldeviation}
\PP(|X- \EE X| \geq \alpha \EE X) \leq 2e^{-\frac{\alpha^2}{3} \EE X}
\end{equation}
In particular, when $X \sim Bi(n,p)$ and $\lambda<\frac{3}{2}np$, then
\begin{equation}\label{chernoff}
\PP(|X-np| \geq \lambda) \leq e^{-\Omega(\lambda^2/(np))}
\end{equation}
\end{lemma}
\begin{lemma} (Chernoff Inequality for large deviation) If $X = \sum_{i=1}^n X_i$, each random variable
$X_i$ has Bernoulli distribution with expectation $p_i$,
and $x \geq 7 ~\EE X$, then
\begin{equation}\label{largedeviation}
\PP(X \geq x) \leq e^{-x}
\end{equation}
\end{lemma}
\textbf{\emph{Proof of Claim \ref{all}}}:
The desired subhypergraph $H''$ is obtained via two rounds of randomization. In the first round, we find edge-disjoint induced subhypergraphs with large minimum degrees which guarantees
the existence of perfect fractional matchings. In the second round, we construct $H''$ from these fractional matchings.

As a preparation toward the first round, $R$ is obtained by choosing
every vertex randomly and independently with probability
$p=|V'|^{-0.9}=n^{-0.9}$. Then $|R|$ is a binomial random variable with expectation $n^{0.1}$. By inequality \eqref{chernoff},
$|R|\sim n^{0.1}$ with  probability $1-e^{-\Omega(n^{0.1})}$.

 Fix a subset $D \subseteq V'$ of size $d$ and let $\DEG_D$ be the number of edges $f\in H'$ such that $D\subset f$ and $f\setminus D\subseteq R$,
 which is the number of edges $e$ in the link graph $H[D]$ with all of its vertices in the random set $R$.
 Therefore $\DEG_D= \sum_{e \in H[D]} X_e$, where $X_e=1$ if $e$ is in $R$ and $0$ otherwise.
We have
\begin{align*}
\EE(\DEG_D) &= \deg_{H'}(D) \times (n^{-0.9})^{k-d} \ge (c+\alpha/2)\binom{n-d}{k-d} n^{-0.9(k-d)} \nonumber \\
&\geq (c+\alpha/3)\binom{|R|-d}{k-d} = \Omega(n^{0.1(k-d)})
\end{align*}
For two distinct intersecting edges $e_i, e_j$ with $|e_i \cap e_j|=l$ for $1 \le l \le k-d-1$, the probability that both of them are in $R$ is
$$\PP(X_{e_i}=X_{e_j}=1)=p^{2(k-d)-l}$$
For fixed $l$, there are at most $\binom{n-d}{k-d}$ choices for $e_i$ in the link graph $H[D]$, $\binom{k-d}{l}$ ways to choose the intersection $L=e_i \cap e_j$ of size $l$,
and $\binom{(n-d)-(k-d)}{k-d-l}$ options for $e_j \backslash L$. Therefore,
\begin{align*}
\Delta &= \sum_{e_i \cap e_j \neq \emptyset} \PP(X_{e_i}=X_{e_j}=1) \leq \sum_{l=1}^{k-d-1} p^{2(k-d)-l} \binom{n-d}{k-d} \binom{k-d}{l} \binom{n-k}{k-d-l}\nonumber \\
&\leq \sum_{l=1}^{k-d-l} p^{2(k-d)-l} O(n^{2(k-d)-l}) = O(n^{0.1(2(k-d)-1)})
\end{align*}
By Janson's inequality (see Theorem 8.7.2 in \cite{AS}),
$$\PP(\DEG_D \le (1-\alpha/12)\EE(\DEG_D)) \le e^{-\Omega((\EE X)^2/\Delta)} \sim e^{-\Omega(n^{0.1})} $$

Therefore by the union bound, with probability $1-n^de^{-\Omega(n^{0.1})}$, for all subsets $D \subseteq V'$ of size $d$, we have
$$ \DEG_D>(1-\alpha/12)\EE(\DEG_D))\ge(c+\alpha/4)\binom{|R|-d}{k-d}.$$

Take $n^{1.1}$ independent copies of $R$ and denote them by $R^i$, $1 \le i \le n^{1.1}$, and the corresponding random variables by $\DEG_D^{(i)}$, where $D \subseteq V'$, $|D|=d$, and $i=1,\dots,n^{1.1}$. Since $|R_i| \sim n^{0.1}$ with probability $1-e^{-\Omega(n^{0.1})}$ for each $i$, the union bound
ensures that $|R_i| \sim n^{0.1}$ for every $i=1,\cdots,n^{1.1}$ with probability $1-o(1)$.
  Now for a subset of vertices $S\subseteq V'$, define the random variable
$$Y_S=|\{i: S\subseteq R^i\}|.$$

Note that the random variables $Y_S$ have binomial distributions $Bi(n^{1.1},n^{-0.9|S|})$ with expectations $n^{1.1-0.9|S|}$.
In particular, for every vertex $v \in V'$, $Y_{\{v\}} \sim Bi(n^{1.1},n^{-0.9})$ and $\EE Y_{\{v\}} = n^{0.2}$.
Hence, by inequality \eqref{chernoff}, taking $\lambda=n^{0.15}$,
$$\PP (|Y_{\{v\}}-n^{0.2}|>n^{0.15}) \leq e^{-\Omega((n^{0.15})^2/n^{0.2})} = e^{-\Omega(n^{0.1})}$$
Therefore a.a.s $|Y_{\{v\}}-n^{0.2}| \leq n^{0.15}$ for every vertex $v \in V'$.

Further, let
$$Z_2=\Big|\left\{\{u,v\}\in \binom {V'}2: Y_{\{u,v\}}\ge3\right\}\Big|.$$ Then
$$\EE Z_2<n^2(n^{1.1})^3(n^{-0.9})^6=n^{-0.1}.$$

Therefore by Markov's inequality,
$$\PP (Z_2 =0) = 1-\PP(Z_2 \geq 1) \geq 1-\EE Z_2> 1-n^{-0.1}$$
This implies that a.a.s every pair of vertices $\{u,v\}$ is contained in at most two subhypergraphs $R^i$.

Finally, for $k\ge3$, let
$$Z_k=\Big|\left\{S\in \binom {V'}k: Y_S\ge2\right\}\Big|.$$ Then,
$$\EE Z_k<n^k(n^{1.1})^2(n^{-0.9})^{2k}=n^{k+2.2-1.8k}\le n^{-0.2}$$
Similarly,
$$\PP (Z_k =0)> 1-n^{-0.2}$$

The latter implies that a.a.s. the induced subhypergraphs $H[R^i]$,  $i=1,\dots,n^{1.1}$, are pairwise edge-disjoint.
Summarizing, we can choose the sets $R^i$, $1 \leq i \leq n^{1.1}$ in such a way that
\begin{enumerate}
\item[(i)] for every $v\in V'$, $Y_{\{v\}} \sim n^{0.2}$,
\item[(ii)] every pair $\{u,v\}\subset V'$  is contained in at most \emph{two} sets $R^i$,
\item[(iii)] every edge $e\in H$  is contained in at most \emph{one} set $R^i$,
\item[(iv)] for all  $i=1,\dots,n^{1.1}$, we have $|R^i|\sim n^{0.1}$, and
\item[(v)] for all  $i=1,\dots,n^{1.1}$ and all $D \subseteq V'$, $|D|=d$, we have $\DEG^{(i)}_D\ge(c+\alpha/4)\binom{|R^i|-d}{k-d}.$
\end{enumerate}
Let us fix a sequence $R^i$, $1 \leq i \leq n^{1.1}$, satisfying (i)-(v) above.

Our assumption that $f_{d}(k,n)\sim c^*\binom{n-d}{k-d}$ holds for all sufficiently large values of $n$, in particular with $n$ replaced by $|R^i| \sim n^{0.1}$. Thus, we have
 $$f_{d}(k,|R^i|)\sim c^*\binom{|R^i|-d}{k-d},$$
and, by condition (v) above, we conclude that
$$\delta_d(H[R^i])\ge(c+\alpha/4)\binom{|R^i|-d}{k-d}>f_{d}(k,|R^i|).$$
 Consequently, by the definition of $f_d$, there exists a fractional perfect matchings $w^i$ in every subhypergraph $H[R^i]$, $i=1,\dots, n^{1.1}$.

Now comes the second round of randomization. Let $H^*=\bigcup_{i}H[R^i]$. We select a generalized binomial subhypergraph $H''$ of $H^*$ by independently choosing each edge $e$ with probability $w^{i_e}(e)$, where $i_e$ is the index $i$ such that $e\in H[R^i]$. Recall that property (iii) ensures that every edge is contained in at most
one hypergraph $R^i$, which guarantees the uniqueness of $i_e$.
We are going to verify our claim by showing $\deg_{H''}(v) \sim n^{0.2}$ for any vertex $v$, while $\Delta_2(H'') \le n^{0.1}$.

 Let $I_v=\{i: v\in R^i\}$ and recall that $|I_v|=Y_{\{v\}}\sim n^{0.2}$ by (i). For every $v\in V'$ the set $E_v$ of edges $e\in H^*$ containing $v$ can be partitioned into $|I_v|$ parts $E_{v}^i=\{e\in E_v\cap H[R^i]\}$.
Recall that $w^i$ is a perfect matching, and thus
$\sum_{e \in E_v^i} w^i(e) = 1 $.
For every $v\in V'$ the random variable $D_v=\deg_{H''}(v)$ is
equal to $\sum_{i\in I_v} \sum_{e\in E_{v}^i} X_e$, where $X_e$ are independent random variables having Bernoulli distribution with expectation $w^{i_e}(e)$.
Therefore $D_v$ is
 generalized binomial with expectation
$$\EE D_v=\sum_{e\in E_v}w^{i_e}(e)=\sum_{i\in I_v}\bigg(\sum_{e\in E_{v}^i}w^i(e)\bigg)=\sum_{i\in I_v} 1\sim n^{0.2}.$$
Hence by Chernoff's inequality \eqref{smalldeviation},
$$\PP(|D_v-n^{0.2}| \geq \alpha n^{0.2} ) \leq 2e^{-\frac{\varepsilon^2}{3}n^{0.2}}$$
Set $\alpha=n^{-0.05}$, then $|D_v-n^{0.2}| \leq n^{0.15}$ with probability $1-O(e^{-n^{0.1}})$.
Taking a union bound over all the $n$ vertices,
we conclude that a.a.s. for all  $v\in V'$ we have $D_v\sim
n^{0.2}$.

Moreover, for all pairs $u,v\in V'$ the random variable
$D_{u,v}=\deg_{H''}(u,v)$  is also generalized binomial with
expectation
$$\EE D_{u,v}=\sum_{e\in E_u\cap E_v}w^{i_e}(e)=\sum_{i\in I_u\cap I_v} \bigg(\sum_{e\in E_{u}^i\cap E_{v}^i}w^i(e) \bigg) \le |I_u \cap I_v| \le 2$$
by (ii). Hence, again by  Chernoff's inequality \eqref{largedeviation} for large deviations, when $n$ is sufficiently large,
$$\PP(D_{u,v} \geq n^{0.1} ) \leq e^{-n^{0.1}}$$
Once again taking the union bound ensures that a.a.s. for every
pair of vertices $u,v\in V'$,
$D_{u,v}\le  n^{0.1}$.
\qed

\section{An application in distributed storage allocation}\label{dstorage}
The following model of distributed storage has been studied in information
theory \cite{dstorage_ldh, dstorage_nr, dstorage_srfs}. A file is split
into multiple chunks,  replicated redundantly and stored
in a distributed storage system with $n$ nodes. Suppose the amount of
data to be stored in each node $i$ is equal to $x_i$, where the
size of the whole file is normalized to $1$. In reality, because
there is limited storage space or transmission bandwidth, we require
that the total amount of data stored does not exceed a given budget $T$,
i.e. $x_1 + \cdots + x_n \leq T$. At the time of retrieval, we attempt
to recover the whole file by accessing only the data stored in a subset $R$ of $r$ nodes which is chosen uniformly at random. It is known that there always exists a coding
scheme such that
we can recover the file whenever the total amount of data accessed is at least $1$.
Our goal is to find an optimal allocation $(x_1, \cdots, x_n)$ in
order to maximize the probability of successful recovery.
This problem can be reformulated as follows.

\begin{question}\rm

For a sequence of nonnegative numbers $(x_1, \cdots, x_n)$, let
$$\Phi(x_1, \cdots, x_n)= \Big|\big\{S \subseteq [n], |S|=r \mbox{ such that } \sum_{i \in S} x_i \geq 1\big\}\Big|.$$
Then the probability of successful recovery of the file equals
$$\frac{\Phi(x_1, \cdots, x_n)}{\binom nr}.$$
Given integers $n \geq r \geq 1$ and a real number $T>0$, determine
$$F^T(r,n) =  \max_{
\sum x_i=T,~x_i \geq 0~\forall i} \Phi(x_1, \cdots, x_n).$$
and find an allocation optimizing $F^T(r,n)$.
\end{question}

In this section, we always assume that $T$ is integer-valued in order to avoid any rounding issues. If the total budget $T$ is at least $ n/r$ then, by setting all $x_i=T/n \geq 1/r$ for all $i$, we can recover the original file
from any subset of size $r$. So, $F^T(r,n)=\binom nr$ for $T \geq n/r$. For $T<n/r$, let $w(i) = x_i$ be a weight function from $V=[n]$ to $\mathbb{R}$.
Then by the definition of the threshold $r$-uniform
hypergraph $H_w^1$ from Section \ref{proofexact},
the edges of $H_w^1$ correspond to the $r$-subsets $S$
such that $\sum_{i \in S} x_i \geq 1$. Thus, it is easy to see that the fractional matching number of $H_w^1$ satisfies
$$\nu^*(H_w^1) = \tau^*(H_w^1) \leq \sum_{i=1}^n w(i) = \sum_{i=1}^n x_i \leq T$$
while $$\Phi(x_1, \cdots, x_n)=|H_w^1|.$$
Therefore, $F^T(r,n)$ is the maximum number of edges in
an $r$-uniform hypergraph on $n$ vertices with fractional matching number at most $T$. As such $F^T(r,n)$ differs from $f_0^T(r,n)$ only in that the latter has the strict inequality $\nu^*(H)<T$ in its definition. But, of course, we have $f_0^T(r,n)\le F^T(r,n)\le f_0^{T+1}(r,n)$, and so $F^T(r,n)\sim f_0^T(r,n)$ as $n\to \infty$.

Hence, Question 5.1 is asymptotically equivalent to the fractional Erd\H os Conjecture \ref{Efrac}.
 As mentioned in
the introduction, it follows from the Erd\H{o}s-Gallai theorem \cite{EG}  that
$$F^T(2,n)\sim f_0^T(2,n) \sim m_0^T(2,n)\sim\max \left\{\binom{2T}{2} ,\binom{n}{2}-\binom{n-T}{2}\right\}.$$
An easy calculation shows that the above maximum equals the first term if $\frac{2}{5}n \le T \leq \frac{1}{2} n$, and the corresponding optimal graph is a clique of size $2T$. This means that, asymptotically, an optimal allocation is $x_1=\cdots = x_{2T}=1/2$ and
$x_{2T+1}=\cdots=x_n=0$. On the other hand, if $T < \frac{2}{5}n$, an optimal allocation is $x_1=\cdots=x_{T}=1$ and $x_{T+1}=\cdots=x_n=0$.

For general $r \geq 3$, if Conjecture \ref{Efrac} is true, then
$$F^T(r,n)\sim\max \left\{ \binom{rT}{r}, \binom{n}{r}-\binom{n-T}{r} \right\}.$$
The bounds are achieved when $H$ is a
clique or a complement of clique. A corresponding (asymptotically) optimal storage allocation is $x_1=\cdots=x_{rT}=1/r, x_{rT+1}=\cdots=x_n=0$ or $x_1=\cdots=x_{T}=1, x_{T+1}=\cdots=x_n=0$, respectively.
Corollary \ref{t95} and Remark \ref{r21} assert that for $r=3$ and
$T<0.277~n$, as well as for $r=4$ and $T<0.217~n$, the latter is an optimal allocation.
Moreover, if Samuels' conjecture \ref{c92} holds for all the remaining $r \geq 5$, then $x_1=\cdots=x_{T}=1, x_{T+1}=\cdots=x_n=0$ is always an asymptotic optimal
allocation whenever $T<n/(r+1)$. Erd\H{o}s \cite{E} proved Conjecture \ref{e}  for all $T< n/(2r^3)$. Recently, the authors of \cite{hls}
extended the range for which this conjecture holds to $T=O(n/r^2)$. Therefore, in this range,  $F^T(r,n)$ is achieved by
the complement of a clique and an optimal allocation is also known to be
$x_1=\cdots=x_{T}=1, x_{T+1}=\cdots=x_n=0$.

\section{Concluding Remarks}\label{concluding}
\begin{list}{\labelitemi}{\leftmargin=1em}
\item We have studied sufficient conditions on the
minimum $d$-degree which guarantee that a uniform hypergraph has
a perfect matching
or perfect fractional matching. We proved that if $f_d(k,n) \sim c^* \binom{n}{k}$, then $m_d(k,n) \sim \max \{c^*, 1/2 \} \binom{n}{k}$.
Therefore in order to determine the asymptotic behavior of the
minimum $d$-degree ensuring existence of a perfect matching,
we can instead
study the presumably easier question for fractional matchings. Using this approach we showed, in particular, that $m_1(5,n)\sim\left(1-\tfrac{4^4}{5^4}\right)\binom{n-1}4$.

\item An intriguing problem which remains open is the conjecture by Erd\H{o}s which states that
the maximum number of edges in a $k$-uniform hypergraph $H$ on $n$ vertices with matching number smaller than $s$
is exactly \begin{equation*} \max \bigg \{\binom{ks-1}{k},~
\binom{n}k-\binom{n-s+1}{k} \bigg\}. \end{equation*}
The fractional version of Erd\H{o}s conjecture is also very interesting. In its
asymptotic form it says that if $H$ is an $l$-uniform $m$-vertex hypergraph with fractional
matching number $\nu^*(H)=xm$, where $0 \leq x <1/l$, then \begin{equation*} \label{erdos} |H| \leq
(1+o(1)) \max \big\{(lx)^l, 1-(1-x)^l \big \} {m \choose l}. \end{equation*}
In Section \ref{proof34} we showed that the fractional Erd\H{o}s conjecture is related to a probabilistic conjecture of
Samuels. This conjecture, if proved, will provide a solution to
the fractional version of Erd\H{o}s problem
for the range $x \leq \frac{1}{l+1}$. It will also lead to
the asymptotics of $m_d(k,n)$ and $f_d(k,n)$ for arbitrary $k \geq d+1$ and $d \geq 1$.

\item As it turns out, matchings and fractional matchings also have
some interesting applications in information theory.
In particular, the uniform model
of distributed storage allocation considered
in \cite{dstorage_srfs} leads to a question which is asymptotically
equivalent to the fractional version of Erd\H{o}s' problem.
In \cite{dstorage_ldh}, the set of accessed nodes, $R$, is given by taking each node randomly and independently with probability $p$.
It would
be interesting to see if our techniques can
be applied to study this binomial model too.

\end{list}

{\bf Acknowledgments} ~The authors would like to thank Alex Dimakis
for a discussion on the fractional Erd\H{o}s conjecture, and an
anonymous referee for many helpful comments.


\end{document}